\documentclass[12pt,oneside,reqno]{amsart}
\usepackage{mathrsfs}
\usepackage{graphics}
\usepackage{amssymb}
\usepackage{enumerate}
\newcommand{\dif}{\mathrm{d}}

\newcommand{\be}{\begin{eqnarray}}
\newcommand{\ee}{\end{eqnarray}}
\newcommand{\ce}{\begin{eqnarray*}}
\newcommand{\de}{\end{eqnarray*}}
\newtheorem{theorem}{Theorem}[section]
\newtheorem{lemma}[theorem]{Lemma}
\newtheorem{remark}[theorem]{Remark}
\newtheorem{definition}[theorem]{Definition}
\newtheorem{proposition}[theorem]{Proposition}
\newtheorem{Example}[theorem]{Example}
\newtheorem{corollary}[theorem]{Corollary}
\newtheorem{condition}[theorem]{Condition}
\def\e{\varepsilon}
\def\s{\sigma}
\def\t{\theta}
\def\a{\alpha}

\def\b{\beta}
\def\d{\delta}

\def\g{\gamma}
\def\l{\lambda}

\def\[{{\Big[}}
\def\]{{\Big]}}
\def\<{{\langle}}
\def\>{{\rangle}}
\def\({{\Big(}}
\def\){{\Big)}}

\def\no{\nonumber}
\def\bt{\begin{theorem}}
\def\et{\end{theorem}}
\def\bl{\begin{lemma}}
\def\el{\end{lemma}}
\def\br{\begin{remark}}
\def\er{\end{remark}}
\def\bx{\begin{Example}}
\def\ex{\end{Example}}
\def\bd{\begin{definition}}
\def\ed{\end{definition}}
\def\bp{\begin{proposition}}
\def\ep{\end{proposition}}
\def\bc{\begin{corollary}}
\def\ec{\end{corollary}}
\def\bco{\begin{condition}}
\def\eco{\end{condition}}
\def\cA{{\mathcal A}}
\def\cB{{\mathcal B}}

\def\cP{{\mathcal P}}

\def\mE{{\mathbb E}}

\def\mN{{\mathbb N}}

\def\mP{{\mathbb P}}
\def\mQ{{\mathbb Q}}
\def\mR{{\mathbb R}}

\def\mW{{\mathbb W}}

\def\sB{{\mathscr B}}
\def\sC{{\mathscr C}}

\def\sF{{\mathscr F}}

\def\sL{{\mathscr L}}

\def\geq{\geqslant}
\def\leq{\leqslant}
\def\epsilon{\varepsilon}

\begin{document}

\allowdisplaybreaks

\title{A limit theorem of nonlinear filtering for multiscale McKean-Vlasov stochastic systems}

\author{Huijie Qiao and Shengqing Zhu}

\dedicatory{School of Mathematics,
Southeast University,\\
Nanjing, Jiangsu 211189, P.R.China}

\thanks{{\it AMS Subject Classification(2020):} 60G35}

\thanks{{\it Keywords:} Multiscale McKean-Vlasov stochastic systems, an average principle, nonlinear filtering, a limit theorem}

\thanks{This work was partly supported by NSF of China (No.12071071).}

\thanks{Corresponding author: Huijie Qiao, hjqiaogean@seu.edu.cn}

\subjclass{}

\date{}

\begin{abstract}
The work concerns about multiscale McKean-Vlasov stochastic systems. First of all, we prove an average principle for these systems in the $L^2$ sense. Moreover, a convergence rate is presented. Then we define the nonlinear filtering of these systems and establish a limit theorem about nonlinear filtering of them in the $L^2$ sense.
\end{abstract}

\maketitle \rm

\section{Introduction}

In this paper, we study the following McKean-Vlasov stochastic system $(X^\e_{t}, Y^\e_{t})$ on $\mR^{n} \times \mR^{m}$: for any $T>0$
\be\left\{\begin{array}{l}
\dif X^\e_{t}=\frac{1}{\e}b(X^\e_{t})\dif t+\frac{1}{\sqrt{\e}}\s(X^\e_{t}) \dif B_{t},\quad 0\leq t\leq T,\\
X^\e_0=\xi,\\
\dif Y^\e_{t}=h(X^\e_{t}, \sL_{X^\e_{t}}^\mQ) \dif t+\dif W_{t},\quad 0\leq t\leq T, \\
Y^\e_0=0,
\end{array}
\right.
\label{xeye}
\ee
where $\{B_{t}, t \geq 0 \}$,$\{W_{t}, t \geq 0\}$ are $d$- and $m$-dimensional standard Brownian motions defined on a complete filtered probability space $(\Omega, \sF, \{\sF_{t}\}_{t \in[0, T]}, \mQ)$, respectively, and mutually independent. Here these coefficients $b: \mR^{n}\mapsto \mR^{n}$, $\s: \mR^{n}\mapsto \mR^{n \times d}$ and $h:\mR^{n} \times \cP_{2}(\mR^{n}) \mapsto \mR^{m}$ are Borel measurable, $\xi$ is a $\sF_0$-measurable random variable with $\mE|\xi|^2<\infty$ and independent of $B.$ and $W.$, and $\sL_{X^\e_{t}}^\mQ$ denotes the distribution of $X^\e_{t}$ under the probability measure $\mQ$ (See Section \ref{pre} for the definition of $\cP_{2}(\mR^{n})$). $\e$ is a small parameter which characterizes the ratio of timescales between processes $X_{\cdot}^{\e}$ and $Y_{\cdot}^{\e}$. 

The system (\ref{xeye}) is usually called a multiscale or fast-slow McKean-Vlasov stochastic system, where $X^\e, Y^\e$ are fast and slow parts, respectively. These systems are widely used in engineering and science fields. For these systems, the average principle is an effective tool in capturing the main behaviors of the slow components while avoiding the complexities caused by the details of the fast components of them. If $h(X^\e_{t}, \sL_{X^\e_{t}}^\mQ)$ is independent of $\sL_{X^\e_{t}}^\mQ$, that is, the system (\ref{xeye}) does not depend on the distribution, about average principles of the system (\ref{xeye}) and more general systems there have been many related results (cf. \cite{rzk, L, LRSX}). If $h(X^\e_{t}, \sL_{X^\e_{t}}^\mQ)$ is replaced by $\tilde{h}(Y^\e_{t}, \sL_{Y^\e_{t}}^\mQ,X^\e_{t})$, where $\tilde h: \mR^{m} \times \cP_{2}(\mR^{m})\times\mR^{n} \mapsto \mR^{m}$ is Borel measurable, there also have been a lot of average principle results on the system (\ref{xeye}) and more general systems (cf. \cite{bs2, lwx, q0, qw1, qw2, rsx, XLLM}). However, for the system (\ref{xeye}), as far as we know, there are few related average principle results. Therefore, in this paper, we establish an average principle for the system (\ref{xeye}). Concretely speaking, we prove that for $0<\g<1$
\ce
\mE\sup\limits_{t\in[0,T]}|Y^\e_t-\bar{Y}_t|^2\leq C(\e+\e^{1-\g}+\e^\g),
\de
where 
\be
\bar{Y}_t:=\bar{h}t+W_t
\label{bary}
\ee
and $\bar{h}:=\int_{\mR^n}h(x,\nu)\nu(\dif x)$ (See Subsection \ref{averprin} for the definition of $\nu$).

Note that in \cite{ghl}, Gao, Hong and Liu considered the following McKean-Vlasov stochastic system
\be\left\{\begin{array}{l}
\dif X^{\e,\d}_{t}=\frac{1}{\e}b(X^{\e,\d}_{t})\dif t+\frac{1}{\sqrt{\e}}\s(X^{\e,\d}_{t}) \dif B_{t},\quad 0\leq t\leq T,\\
X^{\e,\d}_0=\xi,\\
\dif Y^{\e,\d}_{t}=\hat{h}(Y^{\e,\d}_{t}, \sL_{Y^{\e,\d}_{t}}^\mQ, X^{\e,\d}_{t}, \sL_{X^{\e,\d}_{t}}^\mQ) \dif t+\sqrt \d \dif W_{t},\quad 0\leq t\leq T, \\
Y^\e_0=0,
\end{array}
\right.
\label{xeyein}
\ee
where $\hat h: \mR^{m} \times \cP_{2}(\mR^{m})\times\mR^{n}\times \cP_{2}(\mR^{n})\mapsto \mR^{m}$ is Borel measurable and $\d>0$ is a small parameter. There they proved that
\be
\mE\sup\limits_{t\in[0,T]}|Y^{\e,\d}_{t}-\bar{Y}^0_t|^2\leq C(1+\mE|\xi|^2)(\e^{1/3}+\d),
\label{estiin}
\ee
where $C>0$ is a constant independent of $\e, \d$, $\bar{Y}$ solves the following ordinary differential equation
\ce\left\{\begin{array}{l}
\dif\bar{Y}^0_t=\bar{\hat{h}}(\bar{Y}^0_t,\sL_{\bar{Y}^0_t})\dif t,\\
\bar{Y}^0_0=0,
\end{array}
\right.
\de
and $\bar{\hat{h}}(y,\mu):=\int_{\mR^n}\hat{h}(y,\mu,x,\nu)\nu(\dif x)$. Comparing the system (\ref{xeyein}) with the system (\ref{xeye}), we find that the former is more general. But if $\d=1$, by (\ref{estiin}), we only obtain that
\ce
\lim\limits_{\e\rightarrow0}\mE\sup\limits_{t\in[0,T]}|Y^{\e,1}_{t}-\bar{Y}_t|^2\leq C(1+\mE|\xi|^2).
\de
That is, the result in \cite{ghl} does not cover our result (Theorem \ref{yebary}).
 
Next, we define the nonlinear filtering of $(X^\e_{t}, \sL_{X^\e_{t}}^\mQ)$ with respect to $\{Y^\e_{s}, 0 \leq s \leq t\}$ as follows:
\ce
\pi^\e_{t}(F):=\mE[F(X^\e_{t}, \sL_{X^\e_{t}}^\mQ) \mid \sF_{t}^{Y^\e}], \quad F\in\cB_b(\mR^n\times\cP_2(\mR^n)),
\de 
where $\sF_{t}^{Y^\e}$ is the $\mQ$-augmentation $\s$-algebra generated by $\{Y^\e_{s}, 0 \leq s \leq t\}$. And $\pi^\e_{t}(F)$ satisfies a space-distribution dependent Kushner-Stratonovich equation (cf. \cite{lq, sc1, sc2}). Here we show that $\pi^\e_{t}(F)$ converges to $\bar{\pi}(\bar{F})$ in the $L^2$ sense (See Subsection \ref{limitheo} for the definition of $\bar{\pi}(\bar{F})$). Note that if $h, F$ don't depend on the distribution, our result (Theorem \ref{piebarpi}) is the same to \cite[Theorem 3.1 (ii)]{ps2}. Therefore, our result is more general.

The novelty of this work lies in two folds. The first fold is that we prove an average principle for the system (\ref{xeye}) in the $L^2$ sense. Moreover, a convergence rate is presented, which is important for numerical simulation. The second fold is that we establish a limit theorem about nonlinear filtering of the system (\ref{xeye}) in the $L^2$ sense. Comparing Theorem \ref{piebarpi} with some known results (cf. \cite{qw1, qw2}), we find that our result concerns about the nonlinear filtering of the whole signal process and better characterizes the asymptotic behavior of the whole system.

This paper is arranged as follows. In the next section, we introduce some notations used in the sequel. Then the main results are stated in Section \ref{main}. The proofs of two main theorems are placed in Section \ref{yebaryproo} and \ref{piebarpiproo}, respectively. In Section \ref{exam}, we give an example to explain our results.

The following convention will be used throughout the paper: $C$, with or without indices, will denote different positive constants whose values may change from one place to another.

\section{Preliminarie}\label{pre}

In this section, we introduce some notations used in the sequel.
 
Let $|\cdot|$ and $\|\cdot\|$ be norms of vectors and matrices, respectively. Let $A^{*}$ be the transpose of the matrix $A$. 
 
Let $\cB_b(\mR^{n})$ be the set of bounded Borel measurable functions on $\mR^n$. Let $C(\mR^{n})$ be the space of continuous functions on $\mR^{n}$, and $C_{b}(\mR^{n})$ be the collection of functions in $C(\mR^n)$ which are bounded.

Let $\sB(\mR^{n})$ be the Borel $\s$-field on $\mR^{n}$. Let $\cP(\mR^{n})$ be the space of probability measures defined on $\sB(\mR^{n})$ and $\cP_{2}(\mR^{n})$ be the collection of the probability measures $\mu$ on $\sB(\mR^{n})$ satisfying
$$
\mu(|\cdot|^{2}):=\int_{\mR^{n}}|x|^{2} \mu(\dif x)<\infty.
$$
We put on $\cP_{2}(\mR^{n})$ a topology induced by the following $2$-Wasserstein metric:
$$
\mW_{2}^{2}(\mu_{1}, \mu_{2}):=\inf _{\pi \in \sC(\mu_{1}, \mu_{2})} \int_{\mR^{n} \times \mR^{n}}|x-y|^{2} \pi(\dif x, \dif y), \quad \mu_{1}, \mu_{2} \in \cP_{2}(\mR^{n}),
$$
where $ \sC(\mu_{1}, \mu_{2})$ denotes the set of the probability measures whose marginal distributions are $\mu_{1}, \mu_{2}$, respectively. Moreover, if $\xi,\zeta$ are two random variables with distributions $\sL^{\mQ}_\xi, \sL^{\mQ}_\zeta$ under $\mQ$, respectively,
$$
\mathbb{W}^2_2(\sL^{\mQ}_\xi, \sL^{\mQ}_\zeta)\leq\mE|\xi-\zeta|^2,
$$
where $\mE$ stands for the expectation with respect to $\mQ$.

\section{Main results}\label{main}

In this section, we formulate the main results in this paper.

\subsection{The average principle for multiscale McKean-Vlasov stochastic systems}\label{averprin}

In this subsection, we present an average principle for multiscale McKean-Vlasov stochastic systems.

Assume:
\begin{enumerate}[$(\mathbf{H}^1_{b,\s})$]
\item There is a constant $L_1>0$ such that for $x_{1}, x_{2} \in \mR^{n}$,
\ce
|b(x_{1})-b(x_{2})|^2+\|\s(x_{1})-\s(x_{2})\|^{2}\leq L_1|x_{1}-x_{2}|^{2}.
\de
\end{enumerate}
\begin{enumerate}[$(\mathbf{H}_{h})$]
\item $h$ is bounded, and there is a constant $L_2>0$ such that for any $x_i\in \mR^{n}, \mu_i \in \cP_{2}(\mR^{n}), i=1,2$
$$
|h(x_{1},\mu_{1})-h(x_{2},\mu_{2})|^{2}\leq L_2(|x_{1}-x_{2}|^{2}+\mW_2^{2}(\mu_{1},\mu_{2})).
$$
\end{enumerate}
\begin{enumerate}[$(\mathbf{H}^2_{b, \s})$]
\item
There exists a constant $\b>0$ satisfying $\b>2L_1$ such that for $x_{i}\in\mR^n$, $i=1, 2$, 
\ce
2\<x_{1}-x_{2},b(x_1)-b(x_2)\>+\|\s(x_1)-\s(x_2)\|^{2}\leq -\b|x_{1}-x_{2}|^{2}.
\de
\end{enumerate}

\br
By $(\mathbf{H}^1_{b, \s})$ and $(\mathbf{H}^2_{b, \s})$, it holds that for $x\in\mR^n$, 
\be
2\<x,b(x)\>+\|\s(x)\|^{2}\leq -\a|x|^{2}+C,
\label{bemu}
\ee
where $\a:=\b-2L_1$ and $C>0$ is a constant.
\er

Under $(\mathbf{H}^1_{b,\s})$ $(\mathbf{H}_{h})$, the system (\ref{xeye}) has a unique strong solution $(X^\e,Y^\e)$. Then consider the following stochastic differential equation:
\be\left\{\begin{array}{l}
\dif \hat{X}_{t}=b(\hat X_{t})\dif t+\s(\hat X_{t}) \dif \hat{B}_{t},\quad 0\leq t\leq T,\\
\hat X_0=\hat \xi,
\end{array}
\right.
\label{frozequa}
\ee
where $\hat{B}$ is a $d$-dimensional standard Brownian motion defined on the other complete filtered probability space $(\hat{\Omega},\hat{\sF},\{\hat{\sF}_{t}\}_{t \in[0, T]},\hat{\mQ})$ and $\hat \xi$ is a $\hat \sF_0$-measurable random variable independent of $\hat B.$ with $\sL^{\hat \mQ}_{\hat \xi}=\sL^{\mQ}_{\xi}$. Under $(\mathbf{H}^{1}_{b, \s})$, Eq.(\ref{frozequa}) has a unique strong solution $\hat X^{\hat \xi}$. Moreover, under $(\mathbf{H}^2_{b, \s})$, by \cite[Theorem 3.1]{WangFY}, one could obtain that there exists a unique invariant probability measure $\nu\in\cP_2(\mR^n)$ for Eq.(\ref{frozequa}). So, for $\bar{Y}$ defined in (\ref{bary}), we have the following theorem which is the first main result.

\bt\label{yebary}
Suppose that $(\mathbf{H}^1_{b, \s})$, $(\mathbf{H}_{h})$, $(\mathbf{H}^2_{b, \s})$ hold. Then it holds that for $0<\g<1$
\be
\mE\sup\limits_{t\in[0,T]}|Y^\e_t-\bar{Y}_t|^2\leq C(\e+\e^{1-\g}+\e^\g).
\ee
\et

The proof of the above theorem is placed in Section \ref{yebaryproo}.

\br
Let $\g=\frac{1}{2}$, and we obtain that
$$
\left(\mE\sup\limits_{t\in[0,T]}|Y^\e_t-\bar{Y}_t|^2\right)^{1/2}\leq C\e^{1/4}.
$$
That is, the convergence rate is $\frac{1}{4}$.
\er

\subsection{A limit theorem of nonlinear filtering for multiscale McKean-Vlasov stochastic systems}\label{limitheo}

In this subsection, we define nonlinear filtering for multiscale McKean-Vlasov stochastic systems and average systems, and state a limit theorem.

Set 
\ce
(\Lambda^\e_{t})^{-1}:=\exp \left\{-\int_{0}^{t} h^{i}(X^\e_{s}, \sL_{X^\e_{s}}^\mQ) \dif W^i_{s}-\frac{1}{2} \int_{0}^{t}|h(X^\e_{s}, \sL_{X^\e_{s}}^\mQ)|^{2} \dif s\right\}.
\de
Here and hereafter, we use the convention that repeated indices imply summation. Then by $({\bf H}_{h})$, we know that
$(\Lambda^\e_{\cdot})^{-1}$ is an exponential martingale. Define a probability measure $\tilde{\mQ}^\e$ via 
\ce 
\dif \tilde{\mQ}^\e=(\Lambda^\e_{T})^{-1} \dif \mQ,
\de 
and under the probability measure $\tilde{\mQ}^\e$,
$$
Y^\e_{t}=W_{t}+\int_{0}^{t} h(X^\e_{s}, \sL_{X^\e_{s}}^\mQ) \dif s
$$
is an $(\mathscr{F}_t)$-adapted Brownian motion. 

Set for any $F\in\cB_b(\mR^n\times\cP_2(\mR^n))$
\ce
\pi^\e_{t}(F):=\mE[F(X^\e_{t}, \sL_{X^\e_{t}}^\mQ) \mid \sF_{t}^{Y^\e}],
\de 
where $\sF_{t}^{Y^\e}$ is the $\mQ$-augmentation $\s$-algebra generated by $\{Y^\e_{s}, 0 \leq s \leq t\}$, and the Kallianpur-Striebel formula  gives 
\ce
\pi^\e_{t}(F)=\frac{\tilde{\mP}^\e_{t}(F)}{\tilde{\mP}^\e_{t}(1)},
\de
where 
\ce
\tilde{\mP}^\e_{t}(F):=\tilde{\mE}[F(X^\e_{t}, \sL_{X^\e_{t}}^\mQ) \Lambda^\e_{t} \mid \sF_{t}^{Y^\e}], 
\de
and $\tilde{\mE}$ denotes the expectation under the probability measure $\tilde{\mQ}^\e$. Here $\tilde{\mP}_{t}^{\e}(F)$, $\pi_{t}^{\e}(F)$ are called the unnormalized filtering and the normalized filtering of $(X_{t}^{\e},  \sL^{\mQ}_{X_{t}^{\e}})$  with respect to $\mathscr{F}_{t}^{Y^{\e}}$, respectively. 

Also set for any $F\in\cB_b(\mR^n\times\cP_2(\mR^n))$
\ce 
\bar{\Lambda}_{t}:=\exp\left\{\bar{h}^i\bar{Y}^i_t-\frac{1}{2}|\bar{h}|^{2}t \right\}, \quad \bar{F}:=\int_{\mR^n}F(x,\nu)\nu(\dif x),
\de
\ce
\bar{\mP}_t(\bar{F}):=\tilde{\mE}[\bar{F}\bar{\Lambda}_{t}\mid \sF_{t}^{\bar{Y}}],\quad \bar{\pi}_t(\bar{F}):=\frac{\bar{\mP}_t(\bar{F})}{\bar{\mP}_t(1)}.
\de
Then the relationship between $\pi^\e$ and $\bar{\pi}$ is presented in the following theorem which is the main result in this subsection.

\bt\label{piebarpi}
Assume that $(\mathbf{H}^1_{b, \s})$, $(\mathbf{H}_{h})$, $(\mathbf{H}^2_{b, \s})$ hold. Then we have that for any $t\in[0,T]$
\ce 
\lim \limits_{\e\rightarrow0} \mE\left|\pi^{\e}_{t}(F)-\bar{\pi}_{t}(\bar{F})\right|^2=0, \quad F\in \cA_\eta,
\de
where for $\eta>0$ $\cA_\eta$ is defined as follows:
\ce
\cA_\eta&:=&\bigg\{F\in C(\mR^n\times\cP_2(\mR^n)): \sup\limits_{\e\in(0,1)}\sup\limits_{t\in[0,T]}\mE|F(X_t^\e,\sL^\mQ_{X_t^\e})|^{2+\eta}<\infty, \\
&&\qquad\qquad ~\mbox{and}~ \int_{\mR^n}|F(x,\nu)|^2\nu(\dif x)<\infty \bigg\}.
\de
\et

The proof of the above theorem is placed in Section \ref{piebarpiproo}.

\br
Theorem \ref{piebarpi} can be roughly understood as follows: As $\e\rightarrow0$, $\mE[F(X^\e_{t}, \sL_{X^\e_{t}}^\mQ) \mid \sF_{t}^{Y^\e}]$ converges to $``\mE[\bar{F}|\sF_{t}^{\bar{Y}}]"$ in the $L^2$ sense. This is also reasonable.
\er

\br
It is easy to see that $C_b(\mR^n\times\cP_2(\mR^n))\subset\cA_\eta$.
\er

\br
We mention that if both $h$ and the test function $F$ are independent of $\mu$, Theorem \ref{piebarpi} is \cite[Theorem 3.1 (ii)]{ps2}. Therefore, our result is more general.
\er

\section{Proof of Theorem \ref{yebary}}\label{yebaryproo}

In this section, we prove Theorem \ref{yebary}. And we begin with some key estimates.

\bl
Under $(\mathbf{H}^1_{b, \s})$ and $(\mathbf{H}^2_{b, \s})$, it holds that for any $t\geq 0$
\be
&&\hat{\mE}|\hat X_t^{\hat \xi}|^2\leq \hat\mE|\hat\xi|^2e^{-\a t}+\frac{C}{\a}, \label{xbou}\\
&&\hat{\mE}|\hat X_t^{\hat\xi_1}-\hat X_t^{\hat \xi_2}|^2\leq \hat{\mE}|\hat\xi_1-\hat\xi_2|^{2}e^{-\b t}. \label{x1x2esti}
\ee
\el
\begin{proof}
By applying the It\^o formula to $|\hat X_t^{\hat \xi}|^2e^{\l_1 t}$ with $\l_1=\a$, it holds that
\ce
|\hat X_t^{\hat \xi}|^2e^{\l_1 t}&=&|\hat \xi|^2+\l_1\int_0^te^{\l_1 s}|\hat X_s^{\hat \xi}|^2\dif s+2\int_0^te^{\l_1 s}\<\hat X_s^{\hat \xi},b(\hat X_s^{\hat \xi})\>\dif s\\
&&+2\int_0^te^{\l_1 s}\<\hat X_s^{\hat \xi},\s(\hat X_s^{\hat \xi})\dif \hat{B}_s\>+\int_0^te^{\l_1 s}\|\s(\hat X_s^{\hat \xi})\|^2\dif s.
\de
Taking the expectation on two sides, by (\ref{bemu}), we have that
\ce
\hat{\mE}|\hat X_t^{\hat \xi}|^2e^{\l_1 t}&\leq&\hat{\mE}|\hat \xi|^2+\l_1\int_0^te^{\l_1 s}\hat{\mE}|\hat X_s^{\hat \xi}|^2\dif s+\int_0^te^{\l_1 s}\(-\a\hat{\mE}|\hat X_s^{\hat \xi}|^2+C\)\dif s\\
&\leq&\hat{\mE}|\hat \xi|^2+C\int_0^te^{\l_1 s}\dif s,
\de
and
\ce
\hat{\mE}|\hat X_t^{\hat \xi}|^2\leq \hat{\mE}|\hat \xi|^2e^{-\a t}+\frac{C}{\a},
\de
which completes the proof of (\ref{xbou}).

Next, assume that $\hat X_t^{\hat\xi_1}, \hat X_t^{\hat \xi_2}$ are solutions of Eq.(\ref{frozequa}) with the initial values $\hat\xi_1, \hat\xi_2$, respectively. Then it holds that
\ce
\hat X_t^{\hat\xi_1}-\hat X_t^{\hat \xi_2}=\hat\xi_1-\hat\xi_2+\int_{0}^{t}\(b(\hat X_t^{\hat\xi_1})-b(\hat X_t^{\hat\xi_2})\)\dif s+\int_{0}^{t}\(\s(\hat X_t^{\hat\xi_1})
-\s(\hat X_t^{\hat\xi_2})\)\dif \hat{B}_{s}.
\de
Applying the It\^{o} formula to $|\hat X_t^{\hat\xi_1}-\hat X_t^{\hat \xi_2}|^{2}e^{\l_2 t}$ for $\l_2=\b$ and taking the expectation on two sides, we obtain that
\ce
\hat{\mE}|\hat X_t^{\hat\xi_1}-\hat X_t^{\hat \xi_2}|^{2}e^{\l_2 t}
&=&\hat{\mE}|\hat\xi_1-\hat\xi_2|^{2}+\l_2\hat{\mE}\int_{0}^{t}e^{\l_2 s}|\hat X_s^{\hat\xi_1}-\hat X_s^{\hat \xi_2}|^{2}\dif s\\
&&+2\hat{\mE}\int_{0}^{t}e^{\l_2 s}\<\hat X_s^{\hat\xi_1}-\hat X_s^{\hat \xi_2}, b(\hat X_t^{\hat\xi_1})-b(\hat X_t^{\hat\xi_2})\>
\dif s\\
&&+\hat{\mE}\int_{0}^{t}e^{\l_2 s}\|\s(\hat X_t^{\hat\xi_1})
-\s(\hat X_t^{\hat\xi_2})\|^{2}\dif s\\
&\leq&\hat{\mE}|\hat\xi_1-\hat\xi_2|^{2}+(\l_2-\b)\hat{\mE}\int_{0}^{t}e^{\l_2 s}|\hat X_s^{\hat\xi_1}-\hat X_s^{\hat \xi_2}|^{2}\dif s,
\de
which implies that
$$
\hat{\mE}|\hat X_t^{\hat\xi_1}-\hat X_t^{\hat \xi_2}|^{2}\leq \hat{\mE}|\hat\xi_1-\hat\xi_2|^{2}e^{-\b t}.
$$
The proof is complete.
\end{proof}

\bl
Under $(\mathbf{H}^1_{b, \s})$ and $(\mathbf{H}^2_{b, \s})$, it holds that for any $t\geq0$
\be
\mW_2^2(\sL^{\hat{\mQ}}_{\hat X_t^{\hat \xi}},\nu)\leq 2e^{-\b t}(\hat\mE|\hat\xi|^2+\nu(|\cdot|^2)).
\label{xxnu}
\ee
\el
\begin{proof}
Let $\hat\zeta$ be a $\hat{\sF}_0$-measurable random variable with $\sL^{\hat{\mQ}}_{\hat\zeta}=\nu$. And by the definition of $\nu$, we know that $\sL^{\hat{\mQ}}_{\hat X^{\hat\zeta}_t}=\nu$ for any $t\geq0$. So, (\ref{x1x2esti}) implies that
\ce
\mW_2^2(\sL^{\hat{\mQ}}_{\hat X_t^{\hat \xi}},\nu)=\mW_2^2(\sL^{\hat{\mQ}}_{\hat X_t^{\hat \xi}},\sL^{\hat{\mQ}}_{\hat X^{\hat\zeta}_t})\leq \hat{\mE}|\hat X_t^{\hat \xi}-\hat X^{\hat\zeta}_t|^2\leq \hat{\mE}|\hat\xi-\hat\zeta|^2e^{-\b t}\leq 2e^{-\b t}(\hat\mE|\hat\xi|^2+\nu(|\cdot|^2)).
\de
\end{proof}

\bl
Under $(\mathbf{H}^1_{b, \s})$ and $(\mathbf{H}^2_{b, \s})$, it holds that for any $t\geq0$
\be
|\hat{\mE}h(\hat X_t^{\hat \xi},\nu)-\bar{h}|^2\leq 2L_2 e^{-\b t}(\hat\mE|\hat\xi|^2+\nu(|\cdot|^2)).
\label{hbarh}
\ee
\el
\begin{proof}
From the definition of $\nu$ and (\ref{x1x2esti}), it follows that
\ce
|\hat{\mE}h(\hat X_t^{\hat \xi},\nu)-\bar{h}|^2&=&|\hat{\mE}h(\hat X_t^{\hat \xi},\nu)-\int_{\mR^n}h(x,\nu)\nu(\dif x)|^2=|\hat{\mE}h(\hat X_t^{\hat \xi},\nu)-\int_{\mR^n}\hat{\mE}h(\hat X^{x}_t,\nu)\nu(\dif x)|^2\\
&\leq&\hat{\mE}\int_{\mR^n}|h(\hat X_t^{\hat \xi},\nu)-h(\hat X^{x}_t,\nu)|^2\nu(\dif x)\leq L_2\int_{\mR^n}\hat{\mE}|\hat X_t^{\hat \xi}-\hat X^{x}_t|^2\nu(\dif x)\\
&\leq& L_2 e^{-\b t}\int_{\mR^n}\hat\mE|\hat \xi-x|^2\nu(\dif x)\leq 2L_2 e^{-\b t}(\hat\mE|\hat\xi|^2+\nu(|\cdot|^2)),
\de
which completes the proof.
\end{proof}

{\bf Proof of Theorem \ref{yebary}.}

First of all, by (\ref{xeye}) it holds that for any $t\in[0,T]$
\ce
&&|Y^\e_t-\bar{Y}_t|^2=\left|\int_0^t\(h(X^\e_{s}, \sL_{X^\e_{s}}^\mQ)-\bar{h}\)\dif s\right|^2\\
&\leq&2\left|\int_0^t\(h(X^\e_{s}, \sL_{X^\e_{s}}^\mQ)-h(X^\e_{s}, \nu)\)\dif s\right|^2+2\left|\int_0^t\(h(X^\e_{s}, \nu)-\bar{h}\)\dif s\right|^2\\
&=:&I_1+I_2.
\de
For $I_1$, the H\"older inequality and $({\bf H}_{h})$ imply that
\be
\mE\sup\limits_{t\in[0,T]}I_1&\leq& 2TL_2\int_0^T\mW^2_2(\sL_{X^\e_{s}}^\mQ,\nu)\dif s\leq 2TL_2\int_0^T\mW^2_2(\sL_{\hat X^{\hat\xi}_{\frac{s}{\e}}}^{\hat{\mQ}},\nu)\dif s\no\\
&\overset{(\ref{xxnu})}{\leq}&2TL_2\int_0^T2e^{-\b \frac{s}{\e}}(\hat\mE|\hat\xi|^2+\nu(|\cdot|^2))\dif s\no\\
&\leq&4TL_2\b^{-1}(\hat\mE|\hat\xi|^2+\nu(|\cdot|^2))\e.
\label{i1}
\ee

For $I_2$, we divide it into two parts and obtain that
\be
\mE\sup\limits_{t\in[0,T]}I_2&\leq& 4\mE\sup\limits_{t\in[0,T]}\left|\int_0^{[\frac{t}{\d}]\d}\(h(X^\e_{s}, \nu)-\bar{h}\)\dif s\right|^2+4\mE\sup\limits_{t\in[0,T]}\left|\int_{[\frac{t}{\d}]\d}^t\(h(X^\e_{s}, \nu)-\bar{h}\)\dif s\right|^2\no\\
&=&I_{21}+I_{22},
\label{i2122}
\ee
where $\d$ is a fixed positive number depending on $\e$ and $[\frac{t}{\d}]$ denotes the integer part of $\frac{t}{\d}$. 

Let us deal with $I_{21}$. Some computation implies that
\be
I_{21}&\leq& 4\mE\sup\limits_{t\in[0,T]}\left|\sum\limits_{k=0}^{[\frac{t}{\d}]-1}\int_{k\d}^{(k+1)\d}\(h(X^\e_{s}, \nu)-\bar{h}\)\dif s\right|^2\no\\
&\leq&4[\frac{T}{\d}]\mE\sum\limits_{k=0}^{[\frac{T}{\d}]-1}\left|\int_{k\d}^{(k+1)\d}\(h(X^\e_{s}, \nu)-\bar{h}\)\dif s\right|^2\no\\
&\leq&4(\frac{T}{\d})^2\sup\limits_{0\leq k\leq [\frac{T}{\d}]-1}\mE\left|\int_{k\d}^{(k+1)\d}\(h(X^\e_{s}, \nu)-\bar{h}\)\dif s\right|^2\no\\
&=&4(\frac{T}{\d})^2\e^2\sup\limits_{0\leq k\leq [\frac{T}{\d}]-1}\mE\left|\int_{0}^{\d/\e}\(h(X^\e_{s\e+k\d}, \nu)-\bar{h}\)\dif s\right|^2\no\\
&=&8(\frac{T}{\d})^2\e^2\sup\limits_{0\leq k\leq [\frac{T}{\d}]-1}\int_{0}^{\d/\e}\int_{r}^{\d/\e}\mE\<h(X^\e_{s\e+k\d}, \nu)-\bar{h},h(X^\e_{r\e+k\d}, \nu)-\bar{h}\>\dif s\dif r.
\label{i21ces}
\ee
Set $0\leq r\leq s\leq \d/\e$
\ce
\Psi(s,r):=\mE\<h(X^\e_{s\e+k\d}, \nu)-\bar{h},h(X^\e_{r\e+k\d}, \nu)-\bar{h}\>,
\de
and it holds that
\ce
\Psi(s,r)&=&\mE\left[\mE\left[\<h(X^\e_{s\e+k\d}, \nu)-\bar{h},h(X^\e_{r\e+k\d}, \nu)-\bar{h}\>\bigg{|}\sF_{k\d}\right]\right]\no\\
&=&\mE\left[\mE\left[\<h(X^{\e,x}_{s\e}, \nu)-\bar{h},h(X^{\e,x}_{r\e}, \nu)-\bar{h}\>\right]\bigg{|}_{x=X^{\e}_{k\d}}\right],
\de
where $X^{\e,x}_{\e s}$ satisfies the following equation
\ce
X^{\e,x}_{\e s}&=&x+\frac{1}{\e}\int_{0}^{\e s}b(X^{\e,x}_u)\dif u+\frac{1}{\sqrt{\e}}\int_{0}^{\e s}\s(X^{\e,x}_u)\dif B_u.
\de

Besides, we notice that for $0\leq s\leq \d/\e$
\ce
X^{\e,x}_{\e s}=x+\int_{0}^{s}b(X^{\e,x}_{\e v})\dif v+\int_{0}^{s}\s(X^{\e,x}_{\e v})\dif \tilde{B}_v,
\de
where $\tilde{B}_v=\frac{1}{\sqrt{\e}}B_{\e v}$. Since $\tilde{B}$ is a $d$-dimensional standard Brownian motion, $X^{\e,x}_{\e s}$ and the solution $\hat X^{x}_s$ of Eq.(\ref{frozequa}) have the same distributions for $0\leq s\leq \d/\e$. Thus,
\ce
&&\mE\left[\<h(X^{\e,x}_{s\e}, \nu)-\bar{h},h(X^{\e,x}_{r\e}, \nu)-\bar{h}\>\right]\\
&=&\hat{\mE}\<h(\hat X^{x}_s,\nu)-\bar{h}, h(\hat X^{x}_r,\nu)-\bar{h}\>\\
&=&\hat{\mE}\left[\left<\hat{\mE}\left[h(\hat X^{x}_s,\nu)|\sF_r^{\hat{B}}\right]-\bar{h}, h(\hat X^{x}_r,\nu)-\bar{h}\right>\right]\\
&\leq&\left(\hat{\mE}\left|\hat{\mE}\left[h(\hat X^{\t}_{s-r},\nu)\right]|_{\t=\hat X^{x}_r}-\bar{h}\right|^2\right)^{1/2}\left(\hat{\mE}|h(\hat X^{x}_r,\nu)-\bar{h}|^2\right)^{1/2}.
\de

On one hand, by (\ref{hbarh}) and (\ref{xbou}), it holds that
\ce
\left(\hat{\mE}\left|\hat{\mE}\left[h(\hat X^{\t}_{s-r},\nu)\right]|_{\t=\hat X^{x}_r}-\bar{h}\right|^2\right)^{1/2}&\leq& Ce^{-\b (s-r)/2}\left(\hat{\mE}|\hat X^{x}_r|^2+\nu(|\cdot|^2)\right)^{1/2}\\
&\leq& Ce^{-\b (s-r)/2}(1+|x|^2+\nu(|\cdot|^2))^{1/2}.
\de

On the other hand, $({\bf H}_{h})$ and (\ref{xbou}) imply that
\ce
\left(\hat{\mE}|h(\hat X^{x}_r,\nu)-\bar{h}|^2\right)^{1/2}&=&\left(\hat{\mE}|h(\hat X^{x}_r,\nu)-\int_{\mR^n}h(u,\nu)\nu(\dif u)|^2\right)^{1/2}\\
&\leq&\left(\int_{\mR^n}\hat{\mE}|\hat X^{x}_r-u|^2\nu(\dif u)\right)^{1/2}\\
&\leq&C(1+|x|^2+\nu(|\cdot|^2))^{1/2}.
\de
Combining the two calculations, we have that
\ce
\mE\left[\<h(X^{\e,x}_{s\e+k\d}, \nu)-\bar{h},h(X^{\e,x}_{r\e+k\d}, \nu)-\bar{h}\>\right]\leq Ce^{-\b (s-r)/2}(1+|x|^2+\nu(|\cdot|^2)),
\de
and
\ce
\Psi(s,r)&\leq& Ce^{-\b (s-r)/2}(1+\mE|X^{\e}_{k\d}|^2+\nu(|\cdot|^2))=Ce^{-\b (s-r)/2}(1+\hat{\mE}|\hat X^{\hat\xi}_{k\d/\e}|^2+\nu(|\cdot|^2))\\
&\overset{(\ref{xbou})}{\leq}&Ce^{-\b (s-r)/2}(1+\hat\mE|\hat\xi|^2+\nu(|\cdot|^2)),
\de
which together with (\ref{i21ces}) yields that
\be
I_{21}&\leq& 8(\frac{T}{\d})^2\e^2\sup\limits_{0\leq k\leq [\frac{T}{\d}]-1}\int_{0}^{\d/\e}\int_{r}^{\d/\e}Ce^{-\b (s-r)/2}(1+\hat\mE|\hat\xi|^2+\nu(|\cdot|^2))\dif s\dif r\no\\
&\leq&8T^2C(1+\hat\mE|\hat\xi|^2+\nu(|\cdot|^2))\frac{\e}{\d}.
\label{i21}
\ee

Now we treat $I_{22}$. Based on $({\bf H}_{h})$ and the H\"older inequality, it holds that
\be
I_{22}\leq 4\d\mE\sup\limits_{t\in[0,T]}\int_{[\frac{t}{\d}]\d}^t\left|h(X^\e_{s}, \nu)-\bar{h}\right|^2\dif s\leq 4\d\int_0^T\mE\left|h(X^\e_{s}, \nu)-\bar{h}\right|^2\dif s\leq C\d.
\label{i22}
\ee

Finally, collecting (\ref{i1}) (\ref{i21}) (\ref{i22}), one can conclude that
\ce
\mE\sup\limits_{t\in[0,T]}|Y^\e_t-\bar{Y}_t|^2\leq C(\e+\frac{\e}{\d}+\d).
\de
Taking $\d=\e^\g$ for $0<\g<1$, we obtain the required estimate. The proof is complete.

\section{Proof of Theorem \ref{piebarpi}}\label{piebarpiproo}

In this section, we prove Theorem \ref{piebarpi}. First of all, we make some preparations.

\bl\label{4}
Under the assumption $(\mathbf{H}_{h})$, there exists a constant $C>0$ such that
\ce
\mE|\mP_{t}^{\e}(1)|^{-r}\leq \exp\{(2r^{2}+r+1)CT/2\}, \quad r>1.
\de
\el

Since its proof is similar to that of \cite[Lemma 3.6]{q1}, we omit it.

\bp\label{piebarpip}
Assume that $(\mathbf{H}^1_{b, \s})$, $(\mathbf{H}_{h})$, $(\mathbf{H}^2_{b, \s})$ hold. Then we have that for any $t\in[0,T]$
\ce 
\lim \limits_{\e\rightarrow0} \mE\left|\pi^{\e}_{t}(F)-\bar{\pi}_{t}(\bar{F})\right|^2=0, \quad F\in C_{b}(\mR^n\times\cP_2(\mR^n)).
\de
\ep
\begin{proof}
We divide the proof into two parts. In the first step, we prove the main result. In the second step, we show a key convergence which is used in the first step.

{\bf Step 1.} We prove that $\lim\limits_{\e\rightarrow0}\mE|\pi^\e_t(F)-\bar{\pi}_{t}(\bar{F})|^2=0$.

Note that 
\ce
\mE|\pi^\e_t(F)-\bar{\pi}_t(\bar{F})|^2\leq 2\mE|\pi^\e_t(F)-\pi^\e_t(\bar{F})|^2+2\mE|\pi^\e_t(\bar{F})-\bar{\pi}_t(\bar{F})|^2,
\de
and
\ce
\pi^\e_t(\bar{F})=\bar{F}, \quad \bar{\pi}_t(\bar{F})=\frac{\bar{\mP}_t(\bar{F})}{\bar{\mP}_t(1)}=\frac{\bar{F}\bar{\Lambda}_{t}}{\bar{\Lambda}_{t}}=\bar{F}.
\de
Thus, we only need to prove 
\ce
\lim\limits_{\e\rightarrow0}\mE|\pi^\e_t(F)-\pi^\e_t(\bar{F})|^2=0.
\de

On one hand, by the result in {\bf Step 2}, we know that $|\pi^\e_t(F)-\pi^\e_t(\bar{F})|\overset{\mQ}{\rightarrow}0$ as $\e\rightarrow 0$. On the other hand, $F$ and $\bar F$ are bounded. Thus, the bounded dominated convergence theorem yields that 
\ce
\lim\limits_{\e\rightarrow0}\mE|\pi^\e_t(F)-\pi^\e_t(\bar{F})|^2=0.
\de

{\bf Step 2.} We prove that $|\pi^\e_t(F)-\pi^\e_t(\bar{F})|\overset{\mQ}{\rightarrow}0$ as $\e\rightarrow 0$.

By the H\"older inequality, it holds that
\ce
\mE|\pi^\e_t(F)-\pi^\e_t(\bar{F})|=\tilde{\mE}|\pi^\e_t(F)-\pi^\e_t(\bar{F})|\Lambda^\e_T\leq (\tilde{\mE}|\pi^\e_t(F)-\pi^\e_t(\bar{F})|^{r_1})^{1/r_1}(\tilde{\mE}(\Lambda^\e_T)^{r_2})^{1/r_2},
\de
where $1<r_1<2, r_2>1$, $1/r_1+1/r_2=1$. On one hand, it is not difficult to prove that
$$
\(\tilde{\mE}(\Lambda_{T}^{\e})^{r_2}\)^{1/r_2}\leq \exp\{CT\}.
$$
On the other hand, the Kallianpur-Striebel formula and the H\"older inequality imply that
\ce
(\tilde{\mE}|\pi^\e_t(F)-\pi^\e_t(\bar{F})|^{r_1})^{1/r_1}&=&(\tilde{\mE}|\mP^\e_t(F)-\mP^\e_t(\bar{F})|^{r_1}\mP^\e_t(1)^{-r_1})^{1/r_1}\\
&\leq& (\tilde{\mE}|\mP^\e_t(F)-\mP^\e_t(\bar{F})|^{2})^{1/2}(\tilde{\mE}\mP^\e_t(1)^{-2r_1/(2-r_1)})^{\frac{2-r_1}{2r_1}}\\
&\leq&C(\tilde{\mE}|\mP^\e_t(F)-\mP^\e_t(\bar{F})|^{2})^{1/2},
\de
where the last inequality is based on Lemma \ref{4}.

Next, we are devoted to estimating $\tilde{\mE}|\mP^\e_t(F)-\mP^\e_t(\bar{F})|^{2}$. In order to do this, we take an independent copy $\check{X}^\e$ of $X^\e$, which has the same distribution to that for $X^\e$ and is independent of $(X^\e, W)$. It follows from the Jensen inequality that
\be
&&\tilde{\mE}|\mP^\e_t(F)-\mP^\e_t(\bar{F})|^{2}\no\\
&=&\tilde{\mE}\left[\left|\tilde{\mE}\left[F(X_{t}^{\e},\sL^{\mQ}_{X_{t}^{\e}})\Lambda_{t}^{\e}|\mathscr{F}_{t}^{Y^{\e}}\right]-
\tilde{\mE}\left[\bar{F}\Lambda_{t}^{\e}|\mathscr{F}_{t}^{Y^{\e}}\right]\right|^{2}\right]\no\\
&\leq&\tilde{\mE}\left[\tilde{\mE}\left[\left|F(X_{t}^{\e},\sL^{\mQ}_{X_{t}^{\e}})\Lambda_{t}^{\e}
-\bar{F}\Lambda_{t}^{\e}\right|^{2}\bigg{|}\mathscr{F}_{t}^{Y^{\e}}\right]\right]\no\\
&=&\tilde{\mE}\Bigg[\tilde{\mE}\bigg[(F(X_{t}^{\e},\sL^{\mQ}_{X_{t}^{\e}})-\bar{F})(F(\check{X}_{t}^{\e},\sL^{\mQ}_{X_{t}^{\e}})-\bar{F})\no\\
&&\qquad \times\exp\bigg\{\int_{0}^{t} (h^{i}(X^\e_{s}, \sL_{X^\e_{s}}^\mQ)+h^{i}(\check{X}^\e_{s}, \sL_{X^\e_{s}}^\mQ))\dif (Y^\e_{s})^i\no\\
&&\qquad\qquad\qquad -\frac{1}{2} \int_{0}^{t}(|h(X^\e_{s}, \sL_{X^\e_{s}}^\mQ)|^2+|h(\check{X}^\e_{s}, \sL_{X^\e_{s}}^\mQ)|^{2})\dif s\bigg\}\bigg{|}\mathscr{F}_{t}^{Y^{\e}}\bigg]\Bigg]\no\\
&=&\tilde{\mE}\Bigg[\tilde{\mE}\bigg[(F(X_{t}^{\e},\sL^{\mQ}_{X_{t}^{\e}})-\bar{F})(F(\check{X}_{t}^{\e},\sL^{\mQ}_{X_{t}^{\e}})-\bar{F})\no\\
&&\qquad \times\exp\bigg\{\int_{0}^{t} (h^{i}(X^\e_{s}, \sL_{X^\e_{s}}^\mQ)+h^{i}(\check{X}^\e_{s}, \sL_{X^\e_{s}}^\mQ))\dif (Y^\e_{s})^i\no\\
&&\qquad\qquad\qquad -\frac{1}{2} \int_{0}^{t}|h(X^\e_{s}, \sL_{X^\e_{s}}^\mQ)+h(\check{X}^\e_{s}, \sL_{X^\e_{s}}^\mQ)|^{2} \dif s\bigg\}\no\\
&&\qquad \times\exp\bigg\{\int_{0}^{t}h^i(X^\e_{s}, \sL_{X^\e_{s}}^\mQ)h^i(\check{X}^\e_{s}, \sL_{X^\e_{s}}^\mQ)\dif s\bigg\}\bigg{|}\mathscr{F}_{t}^{X^{\e},\check{X}^\e}\bigg]\Bigg]\no\\
&=&\tilde{\mE}\Bigg[(F(X_{t}^{\e},\sL^{\mQ}_{X_{t}^{\e}})-\bar{F})(F(\check{X}_{t}^{\e},\sL^{\mQ}_{X_{t}^{\e}})-\bar{F})\exp\bigg\{\int_{0}^{t}h^i(X^\e_{s}, \sL_{X^\e_{s}}^\mQ)h^i(\check{X}^\e_{s}, \sL_{X^\e_{s}}^\mQ)\dif s\bigg\}\Bigg]\no\\
&=&\tilde{\mE}\Bigg[(F(X_{t}^{\e},\sL^{\mQ}_{X_{t}^{\e}})-\bar{F})(F(\check{X}_{t}^{\e},\sL^{\mQ}_{X_{t}^{\e}})-\bar{F})\no\\
&&\qquad \times\left(\exp\bigg\{\int_{0}^{t}h^i(X^\e_{s}, \sL_{X^\e_{s}}^\mQ)h^i(\check{X}^\e_{s}, \sL_{X^\e_{s}}^\mQ)\dif s\bigg\}-\exp\bigg\{\int_{0}^{t}\bar{h}^ih^i(\check{X}^\e_{s}, \sL_{X^\e_{s}}^\mQ)\dif s\bigg\}\right)\Bigg]\no\\
&&+\tilde{\mE}\Bigg[(F(X_{t}^{\e},\sL^{\mQ}_{X_{t}^{\e}})-\bar{F})(F(\check{X}_{t}^{\e},\sL^{\mQ}_{X_{t}^{\e}})-\bar{F})\exp\bigg\{\int_{0}^{t}\bar{h}^ih^i(\check{X}^\e_{s}, \sL_{X^\e_{s}}^\mQ)\dif s\bigg\}\Bigg]\no\\
&=:&J_1+J_2.
\label{j1j2}
\ee

For $J_1$, by $(\mathbf{H}_{h})$ it holds that
\ce
J_1\leq 8\|F\|^2_{\infty}e^{\|h\|^2_{\infty}T}\tilde{\mE}\left|\int_{0}^{t}(h^i(X^\e_{s}, \sL_{X^\e_{s}}^\mQ)-\bar{h}^i)h^i(\check{X}^\e_{s}, \sL_{X^\e_{s}}^\mQ)\dif s\right|
\de
where $\|F\|_{\infty}$ denotes the boundedness of $F$ and we use the fact that $|e^u-e^v|\leq (e^u+e^v)|u-v|$ for any $u,v\in\mR$. Since $\check{X}^\e$ and $X^\e$ are independent, the ergodicity of $X^{\e}$ yields that 
\be
\lim\limits_{\e\rightarrow0}J_1=0.
\label{j1}
\ee

For $J_2$, we obtain that for arbitrary small $\d>0$
\ce
J_2&\leq& 4\|F\|^2_{\infty}\tilde{\mE}\Bigg|\exp\bigg\{\int_{0}^{t}\bar{h}^ih^i(\check{X}^\e_{s}, \sL_{X^\e_{s}}^\mQ)\dif s\bigg\}-\exp\bigg\{\int_{0}^{t-\d}\bar{h}^ih^i(\check{X}^\e_{s}, \sL_{X^\e_{s}}^\mQ)\dif s\bigg\}\Bigg|\\
&&+\Bigg|\tilde{\mE}\Bigg[\tilde{\mE}\bigg[(F(X_{t}^{\e},\sL^{\mQ}_{X_{t}^{\e}})-\bar{F})\bigg |\mathscr{F}_{t-\d}^{X^{\e}}\vee\sF_t^{\check{X}^\e}\bigg](F(\check{X}_{t}^{\e},\sL^{\mQ}_{X_{t}^{\e}})-\bar{F})\\
&&\qquad\qquad \times \exp\bigg\{\int_{0}^{t-\d}\bar{h}^ih^i(\check{X}^\e_{s}, \sL_{X^\e_{s}}^\mQ)\dif s\bigg\}\Bigg]\Bigg|\\
&\leq& 4\|F\|^2_{\infty}\tilde{\mE}\Bigg|\exp\bigg\{\int_{0}^{t}\bar{h}^ih^i(\check{X}^\e_{s}, \sL_{X^\e_{s}}^\mQ)\dif s\bigg\}-\exp\bigg\{\int_{0}^{t-\d}\bar{h}^ih^i(\check{X}^\e_{s}, \sL_{X^\e_{s}}^\mQ)\dif s\bigg\}\Bigg|\\
&&+2\|F\|_{\infty}\tilde{\mE}\Bigg[\left|\tilde{\mE}\bigg[(F(X_{t}^{\e},\sL^{\mQ}_{X_{t}^{\e}})-\bar{F})\bigg |\mathscr{F}_{t-\d}^{X^{\e}}\bigg]\right|\exp\bigg\{\int_{0}^{t-\d}\bar{h}^ih^i(\check{X}^\e_{s}, \sL_{X^\e_{s}}^\mQ)\dif s\bigg\}\Bigg]\\
&=:&J_{21}+J_{22}.
\de
By $(\mathbf{H}_{h})$ and the dominated convergence theorem, we know that 
\be
\lim\limits_{\d\rightarrow0}J_{21}=0.
\label{j21}
\ee
For $J_{22}$, the ergodicity of $X^{\e}$ implies that for any $\d>0$
\be
\lim\limits_{\e\rightarrow0}J_{22}=0.
\label{j22}
\ee

Finally, combining (\ref{j1}) (\ref{j21}) (\ref{j22}) with (\ref{j1j2}) and first letting $\e\rightarrow0$ and $\d\rightarrow0$, we conclude that
\ce 
\lim\limits_{\e\rightarrow0}\tilde{\mE}|\mP^\e_t(F)-\mP^\e_t(\bar{F})|^{2}=0,
\de
and furthermore
\ce
\lim\limits_{\e\rightarrow0}\mE|\pi^\e_t(F)-\pi^\e_t(\bar{F})|=0.
\de
Since the $L^1$ convergence implies the convergence in probability, $|\pi^\e_t(F)-\pi^\e_t(\bar{F})|$ tends to $0$ in  probability as $\e\rightarrow0$. The proof is complete.
\end{proof}

Now, it is the position to prove Theorem \ref{piebarpi}.

{\bf Proof of Theorem \ref{piebarpi}.} For any $F\in\cA_\eta$, set for $x\in\mR^n, \mu\in\cP_2(\mR^n)$
$$
F_k(x,\mu):=\phi_k(F(x,\mu)), \quad \bar{F}_k=\int_{\mR^n}F_k(x,\nu)\nu(\dif x), \quad k\in\mN,
$$
where $\phi_k(z)=z, |z|\leq k$ and $\phi_k(z)=k{\rm sign}(z), |z|>k$, and it holds that $|F_k(x,\mu)|\leq k$ and
\be
|F_k(x,\mu)-F(x,\mu)|^{2+\eta/2}&\leq&2^{2+\eta/2}|F(x,\mu)|^{2+\eta/2}I_{\{|F(x,\mu)|>k\}}(x,\mu)\no\\
&\leq& 2^{2+\eta/2}k^{-\eta/2}|F(x,\mu)|^{2+\eta}, \label{fkf}\\
|\bar{F}_k-\bar{F}|^2&\leq&\left|\int_{\mR^n}F_k(x,\nu)\nu(\dif x)-\int_{\mR^n}F(x,\nu)\nu(\dif x)\right|^2\no\\
&\leq&4\int_{\mR^n}|F(x,\mu)|^{2}I_{\{|F(x,\mu)|>k\}}(x)\nu(\dif x).\label{barfkbarf}
\ee
Since $F_k\in C_b(\mR^n\times\cP_2(\mR^n))$, Proposition \ref{piebarpip} implies that
\be
\lim\limits_{\e\rightarrow0}\mE|\pi^\e_t(F_k)-\bar{\pi}_t(\bar{F}_k)|^2=0.
\label{pifkbarpifk}
\ee

Besides, by the Jensen inequality and the H\"older inequality, it holds that
\ce
\mE|\pi^\e_t(F)-\pi^\e_t(F_k)|^2&=&\mE\left|\mE[F(X^\e_{t}, \sL_{X^\e_{t}}^\mQ) \mid \sF_{t}^{Y^\e}]-\mE[F_k(X^\e_{t}, \sL_{X^\e_{t}}^\mQ) \mid \sF_{t}^{Y^\e}]\right|^2\\
&\leq&\mE\left|F(X^\e_{t}, \sL_{X^\e_{t}}^\mQ)-F_k(X^\e_{t}, \sL_{X^\e_{t}}^\mQ) \right|^2\\
&\leq&\left(\mE\left|F(X^\e_{t}, \sL_{X^\e_{t}}^\mQ)-F_k(X^\e_{t}, \sL_{X^\e_{t}}^\mQ) \right|^{2+\eta/2}\right)^{\frac{2}{2+\eta/2}}\\
&\overset{(\ref{fkf})}{\leq}&4k^{-\frac{\eta}{2+\eta/2}}\left(\mE\left|F(X^\e_{t}, \sL_{X^\e_{t}}^\mQ)\right|^{2+\eta}\right)^{\frac{2}{2+\eta/2}}\\
&\leq&4k^{-\frac{\eta}{2+\eta/2}}\left(\sup\limits_{\e\in(0,1)}\mE\left|F(X^\e_{t}, \sL_{X^\e_{t}}^\mQ)\right|^{2+\eta}\right)^{\frac{2}{2+\eta/2}},
\de
and furthermore for $\e\in(0,1)$
\be
\lim\limits_{k\rightarrow\infty}\mE|\pi^\e_t(F)-\pi^\e_t(F_k)|^2=0.
\label{pifpifk}
\ee

Finally, we notice that
\ce
\mE|\pi^\e_t(F)-\bar{\pi}_t(\bar{F})|^2&\leq& 3\mE|\pi^\e_t(F)-\pi^\e_t(F_k)|^2+3\mE|\pi^\e_t(F_k)-\bar{\pi}_t(\bar{F}_k)|^2\\
&&+3\mE|\bar{\pi}_t(\bar{F}_k)-\bar{\pi}_t(\bar{F})|^2.
\de
Based on (\ref{barfkbarf}) (\ref{pifkbarpifk}) (\ref{pifpifk}), it holds that 
\ce
\lim\limits_{\e\rightarrow0}\mE|\pi^\e_t(F)-\bar{\pi}_t(\bar{F})|^2=0.
\de
The proof is complete.

\section{An example}\label{exam}

In this section, we give an example to explain our results.

\bx
For $n=d=m=1$, consider the following multiscale McKean-Vlasov stochastic system $(X^\e_{t}, Y^\e_{t})$ on $\mR \times \mR$:
\be\left\{\begin{array}{l}
\dif X^\e_{t}=-\frac{1}{\e}\frac{X^\e_{t}}{2}\dif t+\frac{1}{\sqrt{\e}}\s \dif B_{t},\quad 0\leq t\leq T,\\
X^\e_0=x_0,\\
\dif Y^\e_{t}=\int_{\mR}\sin|X^\e_{t}+u|\sL_{X^\e_{t}}^\mQ(\dif u)\dif t+\dif W_{t},\quad 0\leq t\leq T, \\
Y^\e_0=0,
\end{array}
\right.
\label{xeyeexam}
\ee
where $x_0\in\mR$ and $\s>0$ are two constants. It is easy to see that $b(x)=-x/2, \s(x)=\s, h(x,\mu)=\int_{\mR}\sin|x+u|\mu(\dif u)$ satisfy $(\mathbf{H}^1_{b, \s})$, $(\mathbf{H}_{h})$, $(\mathbf{H}^2_{b, \s})$ with $L_1=\frac{1}{4}, L_2=1, \b=1$. Note that the invariant measure $\nu(\dif x)$ of $X^1$ is $\frac{1}{\sqrt{\pi}\s}e^{-\frac{x}{\s^2}}\dif x$. Thus, by Theorem \ref{yebary}, we obtain that
\ce
\mE\sup\limits_{t\in[0,T]}|Y^\e_t-\bar{Y}_t|^2\leq C(\e+\e^{1-\g}+\e^\g),
\de
where $\bar{Y}_t=\int_{\mR}\int_{\mR}\sin|x+u|\nu(\dif u)\nu(\dif x)t+W_{t}$.

Next, for any $F\in\cB_b(\mR\times\cP_2(\mR))$ define
\ce
\pi^\e_{t}(F):=\mE[F(X^\e_{t}, \sL_{X^\e_{t}}^\mQ) \mid \sF_{t}^{Y^\e}], 
\de
and
\ce
\bar{\Lambda}_{t}:=\exp\left\{\bar{h}^i\bar{Y}^i_t-\frac{1}{2}|\bar{h}|^{2}t \right\}, \quad \bar{F}=\int_{\mR}F(x,\nu)\nu(\dif x),
\de
\ce
\bar{\mP}_t(\bar{F}):=\tilde{\mE}[\bar{F}\bar{\Lambda}_{t}\mid \sF_{t}^{\bar{Y}}], \quad \bar{\pi}_t(\bar{F}):=\frac{\bar{\mP}_t(\bar{F})}{\bar{\mP}_t(1)},
\de 
and by Theorem \ref{piebarpi} it holds that
\ce 
\lim \limits_{\e\rightarrow0} \mE\left|\pi^{\e}_{t}(F)-\bar{\pi}_{t}(\bar{F})\right|^2=0, \quad F\in \cA_\eta.
\de
\ex

\end{document}